\documentclass[11pt,reqno]{amsart}

\usepackage{amssymb}

\newtheorem{Thm}{Theorem}[section]
\newtheorem{Prop}[Thm]{Proposition}
\newtheorem{Lemma}[Thm]{Lemma}
\newtheorem{Cor}[Thm]{Corollary}

\theoremstyle{definition}
\newtheorem{Rem}[Thm]{Remark}

\numberwithin{equation}{section}

\def\endproof{\hfill$\square$\medskip}

\def\beq{\begin{equation}}
\def\eeq{\end{equation}}

\def\C{\mathbb C}
\def\G{\mathbb G}

\def\P{\mathbb P}

\def\Z{\mathbb Z }

\def\Hom{\mbox{Hom}}
\def\dim{\mbox{dim }}

\def\ad{\mbox{ad}}
\def\adt{\mbox{ad}^{+}}

\def\O{\mathcal O }

\def\g{\mathfrak{g}}

\def\P{{P}}

\def\O{\Omega}

\def\O{{\mathcal O}}

\def\<{\langle\kern-.08cm\langle}
\def\>{\rangle\kern-.08cm\rangle}

\def\a{\alpha}
\def\b{\beta}
\def\e{\varepsilon}

\def\gf{\mathfrak{g}}
\def\lf{\mathfrak{l}}
\def\pf{\mathfrak{p}}

\begin{document}

\title[Skew-symmetric degeneracy loci]
{Nonemptiness of skew-symmetric degeneracy loci}
\author{William Graham}
\address{Department of Mathematics\\ University of Georgia\\
Boyd Graduate Studies Research Center\\Athens, GA 30602}
\email{wag@math.uga.edu}
\thanks{Mathematics Subject Classification 14N05.  
Partially supported by the NSF and the Alfred P. Sloan Foundation}

\begin{abstract}
Let $V$ be a rank $N$ vector bundle on a $d$-dimensional complex
projective scheme $X$; assume that $V$ is equipped with a
skew-symmetric bilinear form with values in a line bundle $L$ and that
$\Lambda^2 V^* \otimes L$ is ample.  Suppose that the maximum rank of
the form at any point of $X$ is $r$, where $r>0$ is even.  The main
result of this paper is that if $d>2(N-r)$, then the locus of points
where the rank of the form is at most $r-2$ is nonempty.  The
analogous result for symmetric degeneracy loci was proved in
\cite{Gra:03}; the proof here is similar.  If the hypothesis of
ampleness is relaxed, we obtain a weaker estimate on the maximum
dimension of $X$ (and give a similar result for the symmetric case).
We give applications to subschemes of skew-symmetric matrices, and to
the stratification of the dual of a Lie algebra by orbit dimension.
\end{abstract}

\maketitle

\section{Introduction}
This paper proves a nonemptiness result for skew-symmetric degeneracy
loci.  The analogous result for symmetric degeneracy loci was proved
in \cite{Gra:03}.  The methods used in this paper are similar; the
proof uses ideas of \cite{IlLa:99}, which are related to work of
Fulton, Lazarsfeld, Tu, Harris, and Sommese (\cite{FuLa:81},
\cite{FuLa:83}, \cite{Laz:84}, \cite{Tu:89}, \cite{HaTu:90},
\cite{Som:78}).  The geometry used in the proof is slightly different
than in \cite{Gra:03}: in place of projective and quadric bundles, we
use Grassmannian bundles of $2$-planes and isotropic $2$-planes.  The
results of this paper do not require the more refined properties of
Gysin maps needed to handle the odd rank case in \cite{Gra:03}.

Before stating the main result we illustrate it with an application.
Let $SS_r(N)$ denote the projectivization of the space of
skew-symmetric $N \times N$ complex matrices of rank at most $r$.
Because the rank of a skew-symmetric matrix is necessarily even, we
will assume $r$ is even.  The codimension of $SS_{r-2}(N)$ in
$SS_r(N)$ is $2(N-r)+1$ (cf. \cite[Ex.~14.4.11]{Ful:84}, or
\cite[2.5]{IlLa:99}), so there exist $2(N-r)$-dimensional closed
subschemes of $SS_r(N)$ not meeting $SS_{r-2}(N)$.  We prove that this
is the largest dimension possible.

\begin{Thm}\label{t.application} Assume that $r>0$ is even. If
$X$ is a closed subscheme of $SS_r(N)$ not meeting $SS_{r-2}(N)$,
then $\dim X \le 2(N-r)$.
\end{Thm}

If $U$ is a vector subspace of the space of $N \times N$ complex
matrices, we say that $U$ has constant rank $r$ if every nonzero
matrix in $U$ has rank $r$.  As a corollary to the above theorem,
we obtain the following linear algebra result.

\begin{Cor}\label{c.linalg} Assume $r$ is even.  If $U$ is a constant
rank vector subspace of the space of skew-symmetric $N \times N$
complex matrices, then $\dim U \le 2(N-r)+1$.
\end{Cor}

If $N = r+1$ then the corollary says that $\dim U \le 3$.  The bound
in this case is achieved, as shown by an example of Westwick
(see \cite[p.~168]{IlLa:99}).

If $V \rightarrow X$ is a vector bundle with a bilinear form with
values in a line bundle $L$, let $X_r$ denote the subscheme of points
in $X$ where the rank of the form is at most $r$.  Theorem
\ref{t.application} is deduced from the following theorem, which
is the main result of the paper.

\begin{Thm}\label{t.main} Let $X$ be a $d$-dimensional complex
projective scheme and let $V \rightarrow X$ be a rank $N$ vector
bundle.  Suppose that $V$ is equipped with a skew-symmetric
bilinear form with values in a line bundle $L$ such that the
restriction of the form to any fiber has rank at most $r$, where
$r>0$ is even.  Assume that $\Lambda^2 V^* \otimes L$ is ample.  
If the form has constant rank $r$, then $d \le 2(N-r)$.  
Equivalently, if $d > 2(N-r)$, then the locus $X_{r-2}$ is
nonempty.
\end{Thm}

Theorem \ref{t.application} is an immediate consequence, as
the trivial bundle $V = \C^N \times SS_r(N) \rightarrow SS_r(N)$
is equipped with a skew-symmetric form with values in 
$L = \O(1)$, and the bundle $\Lambda^2 V^* \otimes L$ is ample.
Observe that as a consequence of this theorem one obtains
a stepwise proof of the nonemptiness of skew-symmetric degeneracy
loci, a result proved by other methods in \cite[Prop.~3.5]{FuLa:83}.

Using the notion of the defect of a map defined in
\cite[p.~25]{GoMa:88}, one can define the defect of a vector bundle
(see Section \ref{ss.defect}).  An ample vector bundle has defect $0$.
We can generalize our main theorem to bundles that are not assumed to
be ample:

\begin{Thm}\label{t.mainweak}  Assume the hypotheses of Theorem
\ref{t.main}, except assume that $\Lambda^2 V^* \otimes L$
has defect $e$ (instead of assuming ampleness).  If the form has
constant rank $r$, then $d \le 2(N-r) + e$.  Equivalently, if
$d > 2(N-r) +e$, then the locus $X_{r-2}$ is nonempty.
\end{Thm}

A version of this generalized result also holds in the symmetric
case; see Remark \ref{r.defectsymmetric}.

One natural example of a skew-symmetric map is the Lie bracket on a Lie
algebra $\gf$.  From this, one can define a skew-symmetric form
$\omega$ on the bundle $\gf \times \P(\gf^*) \rightarrow \P(\gf^*)$,
with values in $\O_{\P(\gf^*)}(1)$.
The orbits of the algebraic group $G$ on $\gf^*$ all have even
dimension.  Let
$$
\gf^*_r = \{ \lambda \in \gf^* \ | \ \dim G \cdot \lambda \leq r \};
$$
then $\gf^*_r$ is a closed conical subset of $\gf^*$,
and the projectivization $\P(\gf^*_r)$ coincides with the
locus $\P(\gf^*)_r$ defined using the form $\omega$.  Applying our main
theorem yields a result (Proposition \ref{p.lie}) about the
stratification of $\g^*$ by orbit dimension.  As a consequence of this
result we obtain the following bound on the size of minimal
dimensional orbits in $\g^*$ (Corollary \ref{c.lie}): if $N = \dim \g$,
and $r$ is the minimal dimension of a nonzero $G$-orbit on $\g^*$,
then
$$
r \le 2 \left \lfloor \frac{2N+1}{6} \right \rfloor.
$$
This bound is achieved in the example of the minimal nonzero coadjoint
orbits of the group $SL_3$.

The contents of the paper are as follows.  Section \ref{s.vb} contains
some preliminary results used in the proof of the main theorem.  In
particular, we define projection from a subspace in a Grassmannian (or
Grassmannian bundle) to a smaller Grassmannian, generalizing the
construction for projective space, and we prove that this map is an
affine linear bundle map (Proposition \ref{p.grassaffine}).  The proof
of the main theorem and its generalized version are given in Section
\ref{s.main}, and the application to Lie algebras is given in Section
\ref{s.lie}.

I would like to thank Robert Varley for encouraging me to consider
bundles that are not necessarily ample.

\medskip

{\em Conventions and notation.} Schemes are of finite type over $\C$;
all algebraic groups are assumed to be linear.  Homology and
cohomology are taken with {\em rational} coefficients, unless otherwise
noted.  A symplectic form
on an $r$-dimensional vector space is a nondegenerate skew-symmetric
form; $r$ must be even.  Similarly, we speak of a symplectic form
on a rank $r$ vector bundle; such a form may take values in a line
bundle.  If $V$ and $W$ are vector bundles, $\Hom(V,W)$ denotes the vector
bundle $V^* \otimes W$, and $\G_s(V)$ the Grassmann bundle of
$s$-dimensional planes in $V$.

\section{Preliminaries} \label{s.vb}

Let $r$ be even, and let $( \ , \ )$ be a symplectic form
on $\C^r$.  We define the conformal symplectic group $GSp_r$ as the set of 
all $g$ in $GL_r$ such that for all $v, w \in \C^r$,
\begin{equation} \label{e.conformal}
(gv, gw) = \tau(g) (v,w),
\end{equation}
where $\tau(g)$ depends only on $g$.  Equivalently,
$g$ must satisfy the condition $g^t M g = \tau(g) M$.  
This group is, up to isomorphism, independent of the choice of symplectic
form.  The reason is that
all symplectic forms on $\C^r$ are equivalent, in the sense that
given a second symplectic form on $\C^r$, there is a linear automorphism
of $\C^r$ such that the pullback of the second symplectic form is
our original form.

As in the case of the orthogonal
groups, \eqref{e.conformal} implies that $(\mbox{det }g)^2 = \tau(g)^{r}$. 
However, in the symplectic case, a stronger result holds.  
We include a proof for lack of a reference.

\begin{Prop} \label{p.gsp}
The group $GSp_r$ is connected, and any $g \in GSp_r$ satisfies
$\mbox{det }g = \tau(g)^{r/2}$.
\end{Prop}

\begin{proof}
For $r = 2$, direct computation shows that $GSp_2 = GL_2$ and
$\mbox{det }g = \tau(g)$.  In general, $GSp_r$ acts transitively
on $\P^{r-1}$.  Let $P$ denote the stabilizer in $GSp_r$ of a point
$x$ of $\P^{r-1}$.  Then
$P$ is a parabolic subgroup of $GSp_r$, with Levi factor isomorphic to
$\C^* \times GSp_{r-2}$.  By induction we may assume that $L$ is
connected, and hence so is $P$.  Since $GSp_r /P \simeq \P^{r-1}$
is also connected, we see that $GSp_r$ is connected.  Since
the function
$$
g \mapsto \frac{\tau(g)^{r/2}}{\mbox{det }g}
$$
has square $1$, its only possible values are $\pm 1$.  Hence the
function is constant on connected components; since
$GSp_r$ is connected, the function is identically $1$.
\end{proof}

%

The (homology) Poincar\'e polynomial of a space $X$ is 
by definition $P_t(X) = \sum
b_i t^i$, where $b_i = \dim H_i(X)$ is the $i$-th Betti number of $X$
(recall that we take homology and cohomology with rational
coefficients).  

\begin{Prop} \label{p.poincare}
Let $W$ be a symplectic vector space of dimension $r=2n$, and let $\G_2(W)$
(resp. $Z$) denote the Grassmanian of (resp. isotropic) $2$-planes in $W$.
Then the odd homology groups of $\G_2(W)$ and $Z$ vanish, and, setting
$q = t^2$, the Poincar\'e polynomials of these spaces are given by
\begin{equation}
\begin{array}{ccl}
P_t(\G_2(W))  & = & 1 + q + 2q^2 + 2q^3 + 
   \cdots + (n-1) q^{r-4} + (n-1) q^{r-3}  \\ 
              &   & +  n q^{r-2}  + (n-1) q^{r-1} + \cdots + q^{2r-4} \\
P_t(Z) & = & 1 + q + 2q^2 + 2q^3 + \cdots + 
   (n-1) q^{r-4} + (n-1) q^{r-3} \\ 
       &   & + (n-1) q^{r-2} + (n-1) q^{r-1} + \cdots + q^{2r-5}.
\end{array}
\end{equation}
(Since these spaces are compact orientable manifolds, the polynomials
are palindromic.)
\end{Prop}

\begin{proof}
If $G$ is a connected reductive group, and $P \supset B$ are a parabolic
subgroup and a Borel subgroup, respectively, then there is a fibration
$G/B \rightarrow G/P$, with fibers $P/B$.  Because the odd cohomology
of these spaces vanishes, the cohomology spectral sequence of the
fibration degenerates at $E_2$, so
$$
P_t(G/P) = \frac{P_t(G/B)}{P_t(P/B)}.
$$ The numerator and denominator of the right hand side can be
calculated in terms of Weyl groups. If $l$ is the rank of the
semisimple part of $\gf = \mbox{Lie }G$, then there exist fundamental
invariants $d_1, \ldots, d_l$ of the Weyl group $W$ of $\gf$, called
exponents, in terms of which we can express the Poincar\'e polynomial
of $G/B$:
$$
P_t(G/B) = \sum_{w \in W} q^{l(w)} = \frac{\prod_i (q^{d_i} - 1)}{(q-1)^l}.
$$ Here $l(w)$ is the length of $w$, the first equality follows from
the Bruhat decomposition, and the second is in \cite{Hum:90}.  Because
$P/B = L/B_L$, where the Levi factor $L$ of $P$ is reductive, and
$B_L$ is a Borel subgroup of $L$, $P_t(P/B)$ is calculated by a
similar formula using the exponents of the Weyl group of $\lf =
\mbox{Lie }L$.  The exponents of all Weyl groups of simple complex Lie
algebras are known; for $\mathfrak{sl}_r$ (type $A_{r-1}$) they are
$2,3, \ldots, r$; for $\mathfrak{sp}_r$ (type $C_{r/2}$, with $r$ even), they are
$2,4,6, \ldots, r$ (see \cite{Hum:90}).  If a group is reductive, then
its Lie algebra has a semisimple part which is a product of simple Lie
algebras, and its set of exponents is the union (with multiplicities)
of the sets of exponents corresponding to those simple Lie algebras.

The Grassmannian $\G_2(W)$ equals $G/P$, where $G = SL_r$ and $\pf =
\mbox{Lie }P$ has Levi factor with semisimple part isomorphic to
$\mathfrak{sl}_2 \times \mathfrak{sl}_{r-2}$.  Similarly, the
isotropic Grassmannian $Z$ equals $G/P$, where $G = Sp_r$ and $\pf$ has
Levi factor with semisimple part isomorphic to $\mathfrak{sl}_2
\times \mathfrak{sp}_{r-4}$.  The desired formulas for $P_t(\G_2(W))$ and
$P_t(Z)$ follow (with a little algebra) from this, using the facts in
the previous paragraph.
\end{proof}

Note that as complex varieties, $\dim \G_2(W) = 2r-4$ and $\dim Z = 2r
- 5$, as follows by the above proposition (or by directly computing
the dimensions of the groups $G$ and $P$ used to realize these spaces as
homogeneous spaces).

\begin{Prop} \label{p.betti}
Let $W \rightarrow X$ be a rank $r$ vector bundle over a
$d$-dimensional scheme $X$.  Assume that $W$ is equipped with a
symplectic form with values in a line bundle $L$.  Let $\G_2(W)$
(resp. $Z$) denote the Grassmannian bundles of (resp. isotropic)
$2$-planes in $W$.  Let $b_i = \dim H_i(X)$ for $i \ge 0$, and
$b_i = 0$ for $i<0$.  Then
\begin{equation} \label{e.betti1}
\begin{array}{ccl}
\dim H_{2d+2r-4}(\G_2(W)) & = & n b_{2d} + (n-1) b_{2d-2} + (n-1) b_{2d-4} + 
                         (n-2) b_{2d-6} \\  &  &  + (n-2) b_{2d-8} + \cdots + 
                                 b_{2d-(4r-6)}
                                + b_{2d-(4r-8)} \\
\dim H_{2d+2r-6}(Z) & = &  \dim H_{2d+2r-4}(\G_2(W)) - b_{2d}.
\end{array}
\end{equation}
Hence if $X$ is complete and irreducible, 
\begin{equation} \label{e.betti2}
\dim H_{2d+2r-4}(\G_2(W)) = \dim H_{2d+2r-6}(Z) +1.
\end{equation}
\end{Prop}

\begin{proof} Since both $\G_2(W)$ and $Z$ are partial flag bundles 
associated to principal bundles for connected groups (which are,
respectively, $GL_r$ and $GSp_r$), the homology of each of these
bundles is isomorphic to the tensor product of the homology of $X$
with the homology of the fiber (\cite{Ler:51}; an argument is also
given in \cite[Prop.~4.4]{Gra:03}).  Combining this with the previous
proposition yields \eqref{e.betti1}; \eqref{e.betti2} follows from
this, since if $X$ is complete and irreducible then $b_{2d} = 1$
\cite[Lemma~19.1.1]{Ful:84}.
\end{proof}

We now discuss affine linear bundles on schemes.  By
an affine linear bundle of rank $n$ we mean a fiber bundle with fibers
isomorphic to $\C^n$, but with structure group equal to the group
$M(n)$ generated by $GL_n$ and translations.  We require that such a
bundle be locally trivial in the Zariski topology.  
This definition differs from but is equivalent to the definition
in \cite{Bry:89};
that paper uses the term affine bundle.  However, because this term
is used in \cite[p.~22]{Ful:84} in a weaker sense
(there the structure group is not required to be $M(n)$), we have used
the term affine linear bundle.

We can identify $M(n)$ modulo translations with $GL_n$, and so we have
a natural surjection $\pi: M(n) \rightarrow GL_n$.  Given a collection
of transition functions for an affine linear bundle $E \rightarrow X$,
composing those functions with $\pi$ yields a collection of transition
functions for an associated vector bundle $V \rightarrow X$
(in \cite{Bry:89}, the term associated vector bundle is used
in a different sense).  The
structure group of $E$ reduces to $GL_n$ if and only if $E$ has a
section; in this case, $E$ and its associated vector bundle $V$ are
isomorphic schemes.  

\begin{Rem} \label{r.complextop}
In the usual (complex) topology, $E$ always has a (continuous) section
\cite[12.2]{Ste:51} and thus $E$ and $V$ are homeomorphic.  Therefore
$E$ and $V$ are both homotopy equivalent to $X$.  However, in general
$E$ and $V$ need not be isomorphic as schemes.  For example, if $G =
SL_2$, and $T$ (resp.~$B$) is the subgroup of diagonal 
(resp.~upper triangular) matrices, then
$G/T$ is an affine variety which is an affine linear bundle over $G/B =
\P^1$, and the associated vector bundle is the cotangent bundle of $\P^1$
(see \cite{Bry:89}), which is not an affine variety.
\end{Rem}

Let
$$
0 \rightarrow K \rightarrow V \stackrel{\rho}{\rightarrow} W 
\rightarrow 0
$$
be an exact sequence of vector bundles on a scheme $X$,
of ranks $N-r$, $N$, and $r$, respectively.  
Let $\G_s(V)_K$ denote the open subscheme of $\G_s(V)$ whose fiber over any $x
\in X$ consists of those points in $\G_s(V_x)$ corresponding to those
subspaces of $V_x$ whose intersection with $K_x$ is $\{ 0 \}$. 
Assume that $\G_s(V)_K$ is nonempty, which is equivalent to saying that $s \le
r$.  Define a map
\begin{equation} \label{e.pi}
\pi: \G_s(V)_K \rightarrow \G_s(W)
\end{equation}
by sending $p \in \G_s(V_x)$ to $\rho(p) \in \G_s(W_x)$.  
We call $\pi$ the projection
from the subbundle $K$.

\begin{Prop} \label{p.grassaffine} 
Let 
\begin{equation} \label{e.grassaffine1}
0 \rightarrow K \rightarrow V \stackrel{\rho}{\rightarrow} W
\end{equation}
be an exact sequence of vector bundles on a scheme $X$,
of ranks $N-r$, $N$, and $r$, respectively.  Let $s \leq r$.  
Let $\nu: \G_s(W) \rightarrow X$ denote the projection, and let
$S \rightarrow \G_s(W)$ denote the tautological rank $s$ subbundle
of $\nu^*W$.  Then
$\pi: \G_s(V)_K \rightarrow \G_s(W)$
has the structure of an affine linear bundle on $\G_s(W)$ 
of rank $s(N-r)$, with associated vector bundle
$\Hom (S, \nu^* K)$.
\end{Prop}

We will first prove this in the special case $s = r$.  Note
that $\G_r(W) = X$.  Also, we can identify $\G_r(V)_K$ with the
closed subscheme $E = E(W,K)$
of $\Hom (W,V)$ defined as follows.  
There is a natural map
$\eta: \Hom (W,V) \rightarrow \Hom (W,W)$.  Let $X_1$ denote
the closed subscheme of $\Hom (W,W)$ corresponding to the
identity section, and let 
\begin{equation} \label{e.E}
E = E(W,K) = \eta^{-1}(X_1).
\end{equation}
Note that the fiber $E_x$ is
given by
$$
E_x = \{f \in \Hom(W_x, V_x) \ | \ \rho \circ f = \mbox{id} \} .
$$ 
The isomorphism
$$
E \rightarrow \G_r(V)_K
$$
takes $f \in E_x$ to the subspace $f(W_x) \in \G_r(V_x)$.  The inverse
map takes $U_x \in (\G_r(V)_K)_x$ to the map
$f = (\rho|_{U_x})^{-1} \in E_x$.

As stated in \cite{Ful:96}, $E$ is an affine bundle over $X$.  The
next lemma (the special case $s=r$ of Proposition \ref{p.grassaffine})
shows that $E$ is in fact an affine linear bundle.  In the case where
\eqref{e.grassaffine1} is the tautological exact sequence of bundles
on the Grassmannian $\G_{N-r}(\C^N)$, a version of this lemma is
\cite[Ex.~4.3]{Bry:89}.

\begin{Lemma} \label{l.vectoraffine}
Let $0 \rightarrow K \rightarrow V \stackrel{\rho}{\rightarrow} W
\rightarrow 0$ be an exact sequence of vector bundles on a scheme $X$,
of ranks $N-r$, $N$, and $r$, respectively, and let $E$ be as defined
in \eqref{e.E}.  Then $E$ has a natural structure of an affine linear
bundle on X, with associated vector bundle $Hom(W,K)$.
\end{Lemma}

\begin{proof} We begin by choosing compatible trivializations of the vector
bundles $K$, $V$, and $W$.  Precisely, cover $X$ by (Zariski) open
sets on which we have trivializing isomorphisms
$$
\b_i: V|_{U_i} \stackrel{\cong}{\rightarrow} \C^N \times U_i
$$
such that, if $\a_i$ is the restriction of $\b_i$ to $K|_{U_i}$, then
$\a_i$ takes $K|_{U_i}$ isomorphically to the subbundle
$(\C^{N-r} \times \{ 0 \}) \times U_i$ of $\C^N \times U_i$.
There is an induced isomorphism
$$
\gamma_i: W|_{U_i} \stackrel{\cong}{\rightarrow} \C^r \times U_i.
$$ We denote the transition functions for these bundles by $\a_{ij}$,
$\b_{ij}$, and $\gamma_{ij}$, so, writing $U_{ij} = U_i \cap U_j$ and
we have $\b_{ij}: U_{ij} \rightarrow GL_{N}$, etc.  (To be
precise, write $V|_{ij} = V|_{U_{ij}}$.  Our convention is
that if $S$ is a scheme and $v$ and
$x$ are $S$-valued points of $\C^N$ and $U_{ij}$, respectively, we
define $\b_{ij}(x)$ by the equation
$$
\b_i|_{V_{ij}} \circ (\b_j|_{V_{ij}})^{-1}(v,x) = (\b_{ij}(x)v, x), 
$$ 
and similarly for $\a_{ij}$, $\gamma_{ij}$.)
Note that for a point $x$ of $U_{ij}$, the matrix $\b_{ij}$ has the
block form
$$
\b_{ij}(x) =
\begin{bmatrix} \a_{ij}(x) & M_{ij}(x) \\
0 & \gamma_{ij}(x)
\end{bmatrix}
$$
for some matrix $M_{ij}(x)$.

Our trivializations induce trivializing isomorphisms
$$
\delta_i: \Hom(W,K)|_{U_i} \stackrel{\cong}{\rightarrow} 
\Hom (\C^r, \C^{N-r}) \times U_i.
$$
Let
$$
\delta_{ij}: U_{ij} \rightarrow GL( \Hom (\C^r, \C^{N-r}) )
$$
denote the corresponding transition functions; if
$f \in \Hom (\C^r, \C^{N-r})$ and $x \in U_{ij}$, then
\begin{equation} \label{e.vectoraffine1}
\delta_{ij}(x)(f) = \a_{ij}(x) \circ f \circ \gamma_{ij}(x)^{-1}.
\end{equation}

We have additional induced trivializing isomorphisms
$$
\zeta_i: \Hom(W,V)|_{U_i} \stackrel{\cong}{\rightarrow} 
\Hom (\C^r, \C^{N-r}) \times U_i.
$$
Write $\bar{\zeta_i}$ for the restriction of $\zeta_i$
to $E|_{U_i}$.

Let $p: \Hom (\C^r, \C^{N}) \rightarrow \Hom (\C^r, \C^{N-r})$
be induced by projection of $\C^N$ onto the first
$N-r$ coordinates.  Define maps
$$
\e_i: E|_{U_i} \rightarrow \Hom (\C^r, \C^{N-r}) \times U_i
$$
by
$$
\e_i = (p \times \mbox{id}) \circ \bar{\zeta_i}.
$$
One can verify that the $\e_i$ are isomorphisms and that the
corresponding transition functions $\e_{ij}$ are as follows.  
If $f \in \Hom (\C^r, \C^{N-r})$ and $x \in U_{ij}$, then
\begin{equation} \label{e.vectoraffine2}
\e_{ij}(x)(f) = \a_{ij}(x) \circ f \circ \gamma_{ij}(x)^{-1} +
M_{ij}(x) \gamma_{ij}(x)^{-1}.
\end{equation}
Therefore the $\e_i$ give $E$ the structure of an affine linear
bundle; comparing \eqref{e.vectoraffine1} and \eqref{e.vectoraffine2},
we see that the associated vector bundle is $\Hom (W,K)$.
\end{proof}

\medskip

{\em Proof of Proposition \ref{p.grassaffine}.}  Let 
$\nu: \G_s(W) \rightarrow X$ denote the projection.  The pullback
by $\nu$ of the exact sequence \eqref{e.grassaffine1} is
an exact sequence of vector bundles on $\G_s(W)$:
$$
0 \rightarrow \nu^* K \rightarrow \nu^* V 
\stackrel{\nu^* \rho}{\longrightarrow} \nu^* W \rightarrow 0.
$$
Let $S \subset \nu^* W$ denote the tautological rank $s$ subbundle,
and $B = (\nu^*)^{-1}(S)$.  The sequence
$$
0 \rightarrow \nu^* K \rightarrow B \rightarrow S \rightarrow 0
$$
is again exact.
Let $\tilde{\phi}$ denote the composition
$$
\G_s(B) \rightarrow \G_s(\nu^* V) \cong \G_s(W) \times_X \G_s(V)
\rightarrow \G_s(V),
$$
where the first map is induced by the bundle inclusion
$B \subset \nu^* V$, and the second map is
the natural projection.  Concretely,
if $q \in \G_s(W_x)$ and $p \in \G_s(B_q)$, then the
inclusion
$$
B_q \subset \nu^* V_q = V_{\nu(q)}
$$
means that $p$ defines a point of $\G_s(V_{\nu(q)})$; that
point is $\tilde{\phi}(p)$.  Note that
$\tilde{\phi}^{-1}(\G_s(V)_K) = \G_s(B)_{\nu^* K}$.  
Let 
$$
\phi = \tilde{\phi}|_{\G_s(B)_{\nu^*K}}: \G_s(B)_{\nu^*K} \rightarrow
\G_s(V)_K.
$$
Let $\phi$ denote the restriction of $\tilde{\phi}$ to
$\G_s(B)_{\nu^*K}$.  

We claim that $\phi$ is an isomorphism of schemes over $\G_s(W)$.  By
Lemma \ref{l.vectoraffine}, this suffices, since $\G_s(B)_{\nu^* K}$
can be identified with $E(S, \nu^* K)$.

To prove the claim, by working locally on $X$, we may assume that the
bundles $K$, $V$ and $W$ are trivial, and thereby reduce to the case
where $X$ is a point.  In this case, $K$, $V$, and $W$ are just vector
spaces.  Then
$$
\begin{array}{ccl}
\G_s(B) & \cong & \{ (p,q) \ | \ p \in \G_s(V), \ 
                                 q \in \G_s(W), \
                                 \tilde{q} := \rho^{-1}(q) \supset p \} \\
        & \cong & \{ (\tilde{q} \supset p ) \ | \
                                 p \in \G_s(V), \ 
                                 \tilde{q} \in \G_{s+N-r}(V) , \ 
                                 \tilde{q} \supset K \} .
\end{array}
$$
The projection $\G_s(B) \rightarrow \G_s(W)$ takes
$(\tilde{q} \supset p )$ to $\rho(\tilde{q})$.  The open subvariety
$\G_s(B)_{\nu^*K}$ consists of $(\tilde{q} \supset p )$ as above satisfying
the additional condition $p \cap K = \{ 0 \}$, i.e., 
$p \in \G_s(V)_K$.  The map $\phi$ takes $(\tilde{q} \supset p )$ to $p$.  
The map $p \mapsto (p+K, p)$ is inverse to $\phi$; hence $\phi$ is
an isomorphism.  Finally, $\phi$ is compatible with the
maps to $\G_s(W)$, since given $(\tilde{q} \supset p )$ as above with
$p \in \G_s(V)_K$, we have $\rho(\tilde{q}) = \rho(p)$.  \endproof


\begin{Cor} \label{c.affinez}
Let $V \rightarrow X$ be a rank $N$ vector bundle with a
skew-symmetric form with values in a line bundle $L$.  Assume that the
form has constant rank $r$.  Let $K$ denote the radical of the form,
and let $W = V/K$; then $W$ inherits an $L$-valued skew-symmetric form
$V$.  Let $Z \subset \G_2(W)$ and $\tilde{Z} \subset \G_2(V)$ denote
the subbundles of isotropic $2$-planes.  Then $\G_2(V) - \tilde{Z}$ is
an affine linear bundle over $\G_2(W) - Z$, of rank $2(N-r)$.
\end{Cor}

\begin{proof} This follows from Proposition
\ref{p.grassaffine}, since under the map
$$
\pi: \G_2(V)_K \rightarrow \G_2(W),
$$
the inverse image of $\G_2(W) - Z$ is
$\G_2(V) - \tilde{Z}$.  
\end{proof}

\section{Proofs of the main results} \label{s.main}
\subsection{Proof of Theorem \ref{t.main}} \label{ss.main}
The proof is parallel to that of the main theorem in \cite{Gra:03}.
We may assume that $X$ is irreducible of dimension $d$.  We assume
that the form is of constant rank $r$, where $r$ is even and positive;
we must show that $d \le 2(N-r)$.  Let $K$ denote the radical of the
form; $K$ is a vector subbundle of $V$, of rank $N-r$.  Let $W = V/K$;
then $W$ is a rank $r$ vector bundle on $X$, equipped with a
nondegenerate skew-symmetric form with values in $L$.  Let $\G_2(V)$
(resp. $\G_2(W)$) denote the Grassmann bundle of $2$-planes in $V$
(resp. $W$), and let $\tilde{Z} \subset \G_2(V)$ (resp. $Z \subset
\G_2(W)$) denote the subbundle of isotropic $2$-planes.  By Corollary
\ref{c.affinez}, $\G_2(V) - \tilde{Z}$ is a rank $2(N-r)$ affine
linear bundle over $\G_2(W) - Z$, so by Remark \ref{r.complextop}, the
homology groups of these spaces are isomorphic.

Let $\pi_1$ and $\pi_2$ denote the projections from $\G_2(V)$ and
$\P(\Lambda^2 V)$ to $X$.  We claim that there exists a section of an
ample line bundle on $\G_2(V)$ vanishing only on $\tilde{Z}$.
Indeed, since $\Lambda^2 V^* \otimes L$ is ample, the line bundle
$\O_{\P(\Lambda^2 V \otimes L^*)}(1)$ is ample.  Under the natural
isomorphism between $\P(\Lambda^2 V)$ and $\P(\Lambda^2 V \otimes L^*)$,
this line bundle corresponds to the line bundle $L' = \O_{\P(\Lambda^2
V)}(1) \otimes \pi_2^* L$, which is therefore an ample line bundle on
$\P(\Lambda^2 V)$.  Let $S \rightarrow \G_2(V)$ be the tautological
rank $2$ subbundle.  Under the Pl\"ucker embedding $\G_2(V)
\hookrightarrow \P(\Lambda^2 V)$, $L'$ pulls back to $L'' = \Lambda^2
S^* \otimes \pi_1^* L$, which is therefore ample.  Our skew-symmetric
form is a section of the bundle $\Lambda^2 V^* \otimes L \rightarrow
X$; this pulls back to a section of $\pi_1^*(\Lambda^2 V^* \otimes L)
\rightarrow \G_2(V)$.  Composing this section with the natural map
$\pi_1^*(\Lambda^2 V^* \otimes L) \rightarrow L''$ yields a section of
$L''$.  This section vanishes at $p \in \G_2(V)$ if and only if $p$ is
isotropic, proving the claim.

Because $\G_2(V)$ is projective, the claim implies that $\G_2(V) -
\tilde{Z}$ is an affine scheme.  Since the dimension of $\G_2(V) -
\tilde{Z}$ is $2(N-2) + d$, and since the homology groups of this
space are isomorphic to those of $\G_2(W) - Z$, we conclude that
\begin{equation} \label{e.main1}
H_j(\G_2(W) - Z) = 0 \ \ \mbox{for} \ j > 2(N-2) + d.
\end{equation}
The map $i: Z \hookrightarrow\G_2(W)$ is a regular embedding
and there exists a tubular neighborhood of $Z$ in $\G_2(W)$
(cf. \cite[Prop.~2.5]{Gra:03}).  Therefore there is
a Gysin sequence 
\begin{equation} \label{e.gysin}
\cdots \rightarrow H_j(\G_2(W) - Z) \rightarrow H_j(\G_2(W))
\stackrel{i^*}{\rightarrow} H_{j-2}(Z) \rightarrow H_{j-1}(\G_2(W) - Z) 
\rightarrow \cdots 
\end{equation}
(see \cite[Section 3]{Gra:03}).
This exact sequence and \eqref{e.main1} imply that 
$i^*:H_j(\G_2(W)) \rightarrow H_{j-2}(Z)$ is injective for
$j > 2(N-2) + d$.  On the other hand, 
$$
\dim H_{2d+2r-4}(\G_2(W)) = \dim H_{2d+2r-6}(Z) +1
$$
(Proposition \ref{p.betti}).  Hence $2d+2r-4 \le 2(N-2) + d$,
so $d \le 2(N-r)$, as desired.  This completes the proof.

\begin{Rem} \label{r.gysin}
The preceding proof did not use the more refined properties of 
Gysin maps needed in \cite{Gra:03} to handle the case
of odd rank symmetric degeneracy loci.  
In fact, an exact sequence similar to \eqref{e.gysin} can be constructed
without reference to Gysin maps, and then the proof will go through.
The analogous construction is given in \cite[Section 4]{Gra:03};
we omit details here. 
\end{Rem}

\subsection{Generalization to vector bundles with arbitrary defect} \label{ss.defect}
Our main result generalizes easily to vector bundles with arbitrary
defect, using a result of Goresky and MacPherson.  We begin with some
definitions.  If $\pi: X \rightarrow Y$ is a morphism of varieties,
where $X$ has pure dimension $d$, finitely decompose $X$ into
subvarieties $V_i$ so that $\pi|_{V_i}$ has constant fiber dimension.
Goresky and MacPherson define the defect $D(\pi)$ of $\pi$ as the
supremum over $i$ of the fiber dimension of $\pi|_{V_i}$ minus the
codimension of $V_i$.  They prove that if $\pi$ is proper and $Y$ is
affine, then $X$ has the homotopy type of a CW complex of real
dimension at most $d+e$, where $e = D(\pi)$.  In particular, $H_i(X) =
0$ for $i>d+e$, and $H_{d+e}(X;\Z)$ is torsion-free.  See
\cite[pp.~25,~152]{GoMa:88}.

If $L$ is a line bundle on a complete scheme $X$, we will say that $L$
has defect $D(L) = e$ if some tensor power $L^{\otimes k}$ is pulled
back to $X$ by a morphism $\pi: X \rightarrow \P^l$ with $D(\pi) = e$.
We say a vector bundle $V$ on $X$ has defect $D(V)$ equal to
$D(\O_{\P(V^*)}(1))$.  This is related to the notion of $n$-ampleness
($L$ is $n$-ample if some tensor power $L^{\otimes k}$ is pulled back
to $X$ by a morphism $\pi: X \rightarrow \P^l$ with fiber dimension at
most $n$, and $V$ is $n$-ample if $\O_{\P(V^*)}(1)$ is; see
\cite[Ex.~12.1.5]{Ful:84}).  In fact, an $n$-ample vector bundle has
defect at most $n$; in particular, an ample bundle on a complete
scheme has defect $0$.

The following lemma is an immediate consequence of the result of
Goresky and MacPherson stated above.  For $n$-ample line bundles, a
similar cohomology vanishing theorem appears as \cite[Theorem
6.2]{Tu:90}, with a proof by Harris.

\begin{Lemma} \label{l.homvanish}  Let $L$ be a line bundle on an
irreducible projective scheme $X$ of dimension $d$, and let $Z$ be the
zero-scheme of a non-zero section of $L$.  Assume that $D(L) = e$.
Then $H_i(X-Z) = 0$ for $i > d + e$, and $H_{d+e}(X-Z;\Z)$ is
torsion-free.  \endproof
\end{Lemma}

Using this lemma, the proof of Theorem \ref{t.mainweak} is
a simple modification of the proof of Theorem \ref{t.main};
it is only necessary to replace the estimate
\eqref{e.main1} by the estimate
\begin{equation} \label{e.mainweak}
H_j(\G_2(W) - Z) = 0 \ \ \mbox{for} \ j > 2(N-2) + d +e.
\end{equation}

\begin{Rem} \label{r.defectsymmetric}  In \cite{Gra:03}, we proved
that if $V \rightarrow X$ is a rank $N$ vector bundle on a
$d$-dimensional projective scheme $X$, with an $L$-valued quadratic
form such that $S^2 V^* \otimes L$ is ample, and the quadratic form
has constant rank $r$, then $d \le N-r$.  If instead of ampleness we
assume that $S^2 V^* \otimes L$ has defect $e$, then using Lemma
\ref{l.homvanish}, the proof in \cite{Gra:03} shows that $d \le N-r+e$.
\end{Rem}

\section{Example: Duals of Lie algebras} \label{s.lie}
In this section we apply our main theorem to the stratification of the
dual of a Lie algebra by orbit dimension.

Let $G$ be an algebraic group and $\gf$ its Lie algebra.  Let $\langle
\ , \ \rangle$ denote the pairing between $\g^*$ and $\g$.  Define a
skew-symmetric form $\omega$ on the trivial bundle $\gf \times
\P(\gf^*) \rightarrow \P(\gf^*)$, with values in $\O_{\P(\gf^*)}(1)$,
as follows.  If $\lambda \in \gf^* - \{ 0 \}$, let $[ \lambda ]$
denote the corresponding point in $\P(\gf^*)$.  Then $[ \lambda ]$ is
a line in $\gf$, with dual space $[ \lambda ]^*$.  If $x,y \in \gf$,
then $\omega_{[ \lambda ]}(x,y)$ is the element of $[\lambda]^*$
satisfying
$$
\omega_{[ \lambda ]}(x,y) (\mu)= \langle \mu , [x,y] \rangle,
$$
for $\mu \in [ \lambda ]$.

\begin{Rem}
The form defined above can be viewed as a special case of the
following construction (cf.~ \cite{Gra:03}).  If $V \rightarrow X$ is a
vector bundle with a bilinear form with values in a vector bundle $W$
on $X$, consider $\rho: \P(W^*) \rightarrow X$.  The vector bundle
$\rho^*V$ has a bilinear form with values in $\O_{\P(W^*)}(1)$.  This
is defined by composing the natural $\rho^* W$-valued bilinear form on
$\rho^* V$ with the projection $\rho^* W \rightarrow S^*$, where $S =
\O_{\P(W^*)}(-1)$ is the tautological subbundle.  The form $\omega$ on
the bundle $\g \times \P(\g^*) \rightarrow \P(\g^*)$ is obtained by
applying this construction to $\g$, viewed as a bundle over a point.
\end{Rem}

The group $G$ acts on $\g^*$ via the coadjoint action.
Define 
$$
\gf^*_r = \{ \lambda \in \gf^* \ | \ \dim G \cdot \lambda \leq r \}.
$$ 
This is a conical closed subset of $\gf^*$, so we obtain
$\P(\gf^*_r) \subset \P(\gf^*)$.  

Any coadjoint orbit has a symplectic form (defined by Kirillov and
Kostant) and hence is even-dimensional.  The form $\omega$ is
closely related to the symplectic form on coadjoint orbits; this is
reflected in the following lemma (in which we equip closed
subschemes with the reduced scheme structures).

\begin{Lemma} \label{l.orbitdim}
The subscheme $\P(\gf^*_r)$ of $\P(\gf^*)$ is equal to the subscheme
$\P(\gf^*)_r$ of points in $\P(\gf^*)$ where the form $\omega$
has rank at most $r$.
\end{Lemma}

\begin{proof} Let $G^{\lambda}$ denote the stabilizer of 
$\lambda$ in $G$, and $\g^{\lambda}$ the Lie algebra of $G^{\lambda}$.
It is enough to show that the radical of $\omega_{[ \lambda] }$ is
$\g^{\lambda}$.  Let $\ad$ denote the adjoint action of $\g$ on
itself, and $\adt$ the dual action on $\g^*$.  By definition of the
dual action, for all $x, y \in \g$,
\begin{equation} \label{e.orbitdim}
\langle \adt(x) \lambda, y \rangle + \langle \lambda, \ad(x) y \rangle = 0.
\end{equation}
If $x \in \g$ is in the radical of $\omega_{[ \lambda] }$,
then $\langle \lambda, \ad(x) y \rangle = 0$ for all $y$ in $\g$;
then \eqref{e.orbitdim} implies that $\adt(x) \lambda = 0$, as
desired.
\end{proof}

Applying our main theorem to the stratification of $\P(\g^*)$ yields
as an immediate consequence a result about conical subvarieties of $\g^*$.  

\begin{Prop} \label{p.lie}
Let $G$ be an $N$-dimensional algebraic group with Lie algebra $\g$.  
Let $r>0$ be even, and let 
$\g^*_r$ denote the subscheme of all points in $\g^*$ whose
$G$-orbits have dimension at most $r$.  If $X$ is a conical closed
subscheme of $\g^*_r$ meeting $\g^*_{r-2}$ only at $0$, then 
$\dim X \leq 2(N-r) + 1$. \endproof
\end{Prop}

As a corollary we obtain the following result about minimal orbits.

\begin{Cor} \label{c.lie}
Let $G$ be an $N$-dimensional algebraic group with Lie algebra $\g$.
Let $r$ be the minimal dimension of
a nonzero $G$-orbit in $\g^*$.  Then
$$
r \le 2 \left \lfloor \frac{2N+1}{6} \right \rfloor.
$$
\end{Cor}

\begin{proof} We may assume $r>0$.
Since $\g^*_r$ is a conical closed subscheme of
$\g^*$ which contains subvarieties of dimension $r$ (namely
orbit closures), we have
$$
r \le \dim \g^*_r \leq 2(N-r) + 1,
$$
so 
$\frac{r}{2} \leq \frac{2N + 1}{6}$.  
Because $\frac{r}{2}$ is an integer, the result follows.
\end{proof}

As noted in the introduction, the bound in the above corollary is achieved for
the minimal coadjoint orbits of the group $SL_3$; in this case
$N=8$ and $r=4$.

\def\cprime{$'$}


\begin{thebibliography}{Hum}

\bibitem[Bry]{Bry:89}
Ranee~Kathryn Brylinski, {\em Limits of weight spaces, {L}usztig's
  {$q$}-analogs, and fiberings of adjoint orbits}, J. Amer. Math. Soc.
  \textbf{2} (1989), no.~3, 517--533.

\bibitem[Ful1]{Ful:84}
William Fulton, {\em Intersection theory}, Springer-Verlag, Berlin, 1984.

\bibitem[Ful2]{Ful:96}
\bysame, {\em Schubert varieties in flag bundles for the classical groups},
  Proceedings of the Hirzebruch 65 Conference on Algebraic Geometry (Ramat Gan,
  1993) (Ramat Gan), Israel Math. Conf. Proc., vol.~9, Bar-Ilan Univ., 1996,
  pp.~241--262.

\bibitem[FL1]{FuLa:81}
William. Fulton and Robert Lazarsfeld, {\em On the connectedness of degeneracy
  loci and special divisors}, Acta Math. \textbf{146} (1981), no.~3-4,
  271--283.

\bibitem[FL2]{FuLa:83}
William Fulton and Robert Lazarsfeld, {\em Positive polynomials for ample
  vector bundles}, Ann. of Math. (2) \textbf{118} (1983), 35--60.

\bibitem[GM]{GoMa:88}
Mark Goresky and Robert MacPherson, {\em Stratified {M}orse theory},
  Springer-Verlag, Berlin, 1988.

\bibitem[Gra]{Gra:03}
William Graham, {\em Nonemptiness of symmetric degeneracy loci},
  \rm{arXiv:math.AG/0305159} (2003).

\bibitem[HT]{HaTu:90}
Joe Harris and Loring~W. Tu, {\em The connectedness of symmetric degeneracy
  loci: odd ranks. {A}ppendix to: ``{T}he connectedness of degeneracy loci''
  [{\it {t}opics in algebra, {p}art 2} ({W}arsaw, 1988), 235--248, {PWN},
  {W}arsaw, 1990; {MR} 93g:14050a] by {T}u}, Topics in algebra, Part 2 (Warsaw,
  1988), Banach Center Publ., vol.~26, PWN, Warsaw, 1990, pp.~249--256.

\bibitem[Hum]{Hum:90}
James~E. Humphreys, {\em Reflection groups and {C}oxeter groups}, Cambridge
  Studies in Advanced Mathematics, vol.~29, Cambridge University Press,
  Cambridge, 1990.

\bibitem[IL]{IlLa:99}
Bo~Ilic and J.~M. Landsberg, {\em On symmetric degeneracy loci, spaces of
  symmetric matrices of constant rank and dual varieties}, Math. Ann.
  \textbf{314} (1999), 159--174.

\bibitem[Laz]{Laz:84}
Robert Lazarsfeld, {\em Some applications of the theory of positive vector
  bundles}, Complete intersections (Acireale, 1983), Lecture Notes in Math.,
  vol. 1092, Springer, Berlin, 1984, pp.~29--61.

\bibitem[Ler]{Ler:51}
Jean Leray, {\em Sur l'homologie des groupes de {L}ie, des espaces homog\`enes
  et des espaces fibr\'es principaux}, Colloque de topologie (espace fibr\'es),
  Bruxelles, 1950, Georges Thone, Li\`ege, 1951, pp.~101--115.

\bibitem[Som]{Som:78}
Andrew~John Sommese, {\em Submanifolds of {A}belian varieties}, Math. Ann.
  \textbf{233} (1978), no.~3, 229--256.

\bibitem[Ste]{Ste:51}
Norman Steenrod, {\em The {T}opology of {F}ibre {B}undles}, Princeton
  Mathematical Series, vol. 14, Princeton University Press, Princeton, N. J.,
  1951.

\bibitem[Tu1]{Tu:89}
Loring~W. Tu, {\em The connectedness of symmetric and skew-symmetric degeneracy
  loci: even ranks}, Trans. Amer. Math. Soc. \textbf{313} (1989), no.~1,
  381--392.

\bibitem[Tu2]{Tu:90}
\bysame, {\em The connectedness of degeneracy loci}, Topics in algebra, Part 2
  (Warsaw, 1988), Banach Center Publ., vol.~26, PWN, Warsaw, 1990,
  pp.~235--248.

\end{thebibliography}
\end{document}